\documentclass[a4paper,12pt]{article}
\usepackage{latexsym}
\usepackage{amsmath}
\usepackage{amssymb}
\usepackage{amscd}
\usepackage{array}
\usepackage{comment}
\usepackage[all]{xy}
\usepackage{amsthm}
\usepackage{amsfonts}
\usepackage{enumerate}
\usepackage{hyperref}
\usepackage{tabu}
\usepackage{multirow}
\usepackage[normalem]{ulem}
\usepackage{pdflscape}
\usepackage{authblk}
\usepackage{pdfpages}
\newtheorem{theorem}{Theorem}[]
\newtheorem{proposition}[theorem]{Proposition}

\theoremstyle{definition}

\newtheorem{remark}[theorem]{Remark}

\newcommand{\Z}{\mathbf Z}

\newcommand{\F}{\mathbf F}

\newcommand{\Gal}{\mathrm{Gal}}

\newcommand{\Hol}{\mathrm{Hol}}

\newcommand{\Sym}{\operatorname{Sym}}

\newcommand{\GL}{\mathrm{GL}}

\newcommand{\End}{\operatorname{End}}

\newcommand{\Aut}{\operatorname{Aut}}

\newcommand{\wL} {{\widetilde{L}}}

\title{An algorithm to determine Hopf Galois structures}
\author[1]{Teresa Crespo}

\author[2]{Marta Salguero}

\affil[1,2]{\small{Departament de Matem\`atiques i Inform\`atica, Universitat de Barcelona (UB), \newline Gran Via de les
Corts Catalanes 585, E-08007 Barcelona, Spain, \newline e-mail: teresa.crespo@ub.edu, msalguga11@alumnes.ub.edu}}

\date{\today}

\begin{document}

\maketitle

\let\thefootnote\relax\footnotetext{Both authors acknowledge support by grant MTM2015-66716-P (MINECO/FEDER, UE).}

\begin{abstract} A Hopf Galois structure on a finite field extension $L/K$ is a pair $(H,\mu)$, where $H$ is a finite cocommutative $K$-Hopf algebra and $\mu$ a Hopf action. In this paper we present an algorithm written in the computational algebra system Magma which gives all Hopf Galois structures on separable field extensions of a given degree and several properties of those. We describe the results obtained for extensions of degree up to 11. Besides, we prove that separable extensions of degree $p^2$, for $p$ an odd prime, have at most one type of Hopf Galois structures.

\end{abstract}

\section{Introduction}
A Hopf Galois structure on a finite extension of fields $L/K$ is a pair $(H,\mu)$, where $H$ is
a finite cocommutative $K$-Hopf algebra  and $\mu$ is a
Hopf action of $H$ on $L$, i.e a $K$-linear map $\mu: H \to
\End_K(L)$ giving $L$ a left $H$-module algebra structure and inducing a $K$-vector space isomorphism $L\otimes_K H\to\End_K(L)$.
Hopf Galois structures were introduced by Chase and Sweedler in \cite{C-S}.
For separable field extensions, Greither and
Pareigis \cite{G-P} give the following group-theoretic
equivalent condition to the existence of a Hopf Galois structure.

\begin{theorem}\label{G-P}
Let $L/K$ be a separable field extension of degree $g$, $\wL$ its Galois closure, $G=\Gal(\wL/K), G'=\Gal(\wL/L)$. Then there is a bijective correspondence
between the set of Hopf Galois structures on $L/K$ and the set of
regular subgroups $N$ of the symmetric group $S_g$ normalized by $\lambda (G)$, where
$\lambda:G \hookrightarrow S_g$ is the monomorphism given by the action of
$G$ on the left cosets $G/G'$.
\end{theorem}

For a given Hopf Galois structure on a separable field extension $L/K$ of degree $g$, we will refer to the isomorphism class of the corresponding group $N$ as the type of the Hopf Galois
structure. The Hopf algebra $H$ corresponding to a regular subgroup $N$ of $S_g$ normalized by $\lambda (G)$ is the Hopf subalgebra $\wL[N]^G$ of the group algebra $\wL[N]$ fixed under the action of $G$, where $G$ acts on $\wL$ by $K$-automorphisms and on $N$ by conjugation through $\lambda$. The Hopf action is induced by $n \mapsto n^{-1}(\overline{1})$, for $n \in N$, where we identify $S_g$ with the group of permutations of $G/G'$ and $\overline{1}$ denotes the class of $1_G$ in $G/G'$. It is known that the Hopf subalgebras of $\wL[N]^G$ are in 1-to-1 correspondence with the subgroups of $N$ stable under the action of $G$ (see e.g. \cite{CRV} Proposition 2.2) and that, given two regular subgroups $N_1, N_2$ of $S_g$ normalized by $\lambda (G)$, the Hopf algebras $\wL[N_1]^G$ and $\wL[N_2]^G$ are isomorphic if and only if the groups $N_1$ and $N_2$ are $G$-isomorphic.

Childs \cite{Ch1} gives an equivalent  condition to the existence of a Hopf Galois structure introducing the holomorph of the regular subgroup $N$ of $S_g$. We state the more precise formulation of this result due to Byott \cite{B} (see also \cite{Ch2} Theorem 7.3).

\begin{theorem}\label{theoB} Let $G$ be a finite group, $G'\subset G$ a subgroup and $\lambda:G\to \Sym(G/G')$ the morphism given by the action of
$G$ on the left cosets $G/G'$.
Let $N$ be a group of
order $[G:G']$ with identity element $e_N$. Then there is a
bijection between
$$
{\cal N}=\{\alpha:N\hookrightarrow \Sym(G/G') \mbox{ such that
}\alpha (N)\mbox{ is regular}\}
$$
and
$$
{\cal G}=\{\beta:G\hookrightarrow \Sym(N) \mbox{ such that }\beta
(G')\mbox{ is the stabilizer of } e_N\}
$$
Under this bijection, if $\alpha\in {\cal N}$ corresponds to
$\beta\in {\cal G}$, then $\alpha(N)$ is normalized by
$\lambda(G)$ if and only if $\beta(G)$ is contained in the
holomorph $\Hol(N)$ of $N$.
\end{theorem}

In Hopf Galois theory one has the following Galois correspondence theorem.

\begin{theorem}[\cite{C-S} Theorem 7.6]\label{esto} Let $(H,\mu)$ be a Hopf Galois structure on the field extension $L/K$.
For a $K$-sub-Hopf algebra $H'$ of $H$ we define
$$
L^{H'}=\{x\in L \mid \mu(h)(x)=\varepsilon(h)\cdot x \mbox{ for all } h\in H'\},
$$
where $\varepsilon$ is the counity of $H$.
Then, $L^{H'}$ is a subfield of $L$, containing $K$, and
$$
\begin{array}{rcl}
{\mathcal F}_H:\{H'\subseteq H \mbox{ sub-Hopf algebra}\}&\longrightarrow&\{\mbox{Fields }E\mid K\subseteq E\subseteq L\}\\
H'&\to &L^{H'}
\end{array}
$$
is injective and inclusion reversing.
\end{theorem}

In \cite{G-P} a class of Hopf Galois structures is identified for which the Galois correspondence is bijective. We shall say that a Hopf Galois structure $(H,\mu)$ on
$L/K$ is an \emph{almost classically Galois structure} if the corresponding regular subgroup $N$ of $S_g$ normalized by $\lambda(G)$ satisfies that its centralizer $Z_{S_g}(N)$ in $S_g$ is contained in $\lambda(G)$.

\begin{theorem}[\cite{G-P} 5.2]
If $(H,\mu)$ is an almost classically Galois Hopf Galois structure on $L/K$, then the map ${\mathcal F}_H$ from the set of $K$-sub-Hopf algebras of $H$ into the set of subfields of $L$ containing $K$ is bijective.
\end{theorem}

In \cite{CRV2} the Hopf Galois character of separable field extensions of degree up to 7 and of some subextensions of its normal closure has been determined. In \cite{CRV} Theorem 3.4, a family of extensions is given with no almost classically Galois structure but with a Hopf Galois structure for which the Galois correspondence is bijective. In \cite{CRV4} a degree 8 non-normal separable extension having two non-isomorphic Hopf Galois structures with isomorphic underlying
Hopf algebras is presented.

In this paper we present an algorithm which determines all Hopf Galois structures of a separable field extension of a given degree $g$ and their corresponding type. Moreover our algorithm distinguishes almost classically Galois structures and decides for the remaining ones if the Galois correspondence is bijective. Finally it classifies the Hopf Galois structures in Hopf algebra isomorphism classes. We detail the results obtained for separable field extensions of composite degree up to 11. In the case of prime degree, we obtain the results already found in \cite{Ch1} theorem 2 and \cite{P} theorem 5.2, namely that if $L/K$ is a separable field extension of prime degree and $\wL$ its Galois closure, then $L/K$ has a Hopf Galois structure if and only if $\Gal(\widetilde{L}/K)$ is solvable and, in this case, the Hopf Galois structure is unique. We note that the case of degree 8 is especially interesting since there are 5 groups of order 8, up to isomorphism. In the case of degree $p^2$, for $p$ an odd prime, we prove that a separable field extension of degree $p^2$ has at most one type of Hopf Galois structures.

\section{Extensions of degree $p^2$, for $p$ an odd prime}

For $p$ prime, there are exactly two groups of order $p^2$, up to isomorphism, the cyclic one $C_{p^2}$ and the direct product of two copies of $C_p$, hence two possible types for a Hopf Galois structure of a field extension of degree $p^2$. We shall prove that the two types do not occur simultaneously, when $p\neq 2$. The case $p=2$ goes differently, see table \ref{d4-ext}. For $p=3$, we give in table \ref{d9-ext} the results obtained by our algorithm. If we write $C_{p^2}$ additively as $\Z/p^2 \Z$, its holomorph is $\Z/p^2 \Z \rtimes (\Z/p^2 \Z)^*$. For $C_p\times C_p$ the automorphism group is isomorphic to $\GL(2,\F_p)$.

\begin{proposition}\label{psqua}
Let $L/K$ be a separable field extension of degree $p^2$, $p$ an odd prime, $\widetilde{L}/K$ its normal closure and $G\simeq \Gal(\widetilde{L}/K)$.
If $L/K$ has a Hopf Galois structure of cyclic type, then it has no structure of type $C_p\times C_p$. Therefore a separable field extension of degree $p^2$, $p$ an odd prime,  has at most one type of Hopf Galois structures either cyclic or $C_p\times C_p$.
\end{proposition}

\noindent {\it Proof.}
 By theorem \ref{theoB}, if $L/K$ has a Hopf Galois structure of type $C_{p^2}$, then $G$ is a transitive subgroup of $\Hol(C_{p^2})$. We shall see that all transitive subgroups of $\Hol(C_{p^2})$ contain an element of order $p^2$. Let us write $\Hol(C_{p^2})$ as $\Z/p^2 \Z \rtimes (\Z/p^2 \Z)^*$ and let $\sigma$ be a generator of $(\Z/p^2 \Z)^*$. The immersion of  $\Hol(C_{p^2})$ in the symmetric group $S_{p^2}$ is given by sending the generator $1$ of $\Z/p^2 \Z$ to the $p^2$-cycle $(1,2,\dots,p^2)$ and $\sigma$ to itself, considered as a permutation. The stabilizer of 1 in the image $H$ of  $\Hol(C_{p^2})$ in $S_{p^2}$ consists in the elements $(1-k^j,\sigma^j)$, where $k=\sigma(1)$ and $1-k^j$ is computed modulo $p^2$. We have $|H|=|\Hol(C_{p^2})|=p^3(p-1)$, hence $H$ has a unique $p$-Sylow subgroup $Syl(H)$ which is isomorphic to the only non-abelian group of order $p^3$ having an element of order $p^2$ (see \cite{C}). Now, a subgroup $H'$ of $H$ is transitive if and only if $[H':Stab_H(1) \cap H']=p^2$. Let $H'$ be a transitive subgroup of $H$. We have then $p^2 \mid |H'|$ and $|H'| \mid p^3(p-1)$, hence $H'$ has a unique $p$-Sylow subgroup $Syl(H')$ which has order $p^3$ or $p^2$. In the first case, $Syl(H')=Syl(H)$ contains an element of order $p^2$. In the second case, $Syl(H')$ is a subgroup of $Syl(H)$ of order $p^2$. The group $Syl(H)$ is isomorphic to the group

 $$G_p:= \left\{ \left(\begin{array}{cc} 1+pm & b \\ 0 & 1 \end{array} \right) : m, b \in \Z/p^2 \Z \right\},$$

 \noindent where $m$ actually only matters modulo $p$. The group $G_p$ has $p^3-p^2$ elements of order $p^2$, those with $b\not \equiv 0 \pmod{p}$, hence $p$ cyclic subgroups of order $p^2$ and $p^2-1$ elements of order $p$, those nontrivial with $b \equiv 0 \pmod{p}$, hence one noncyclic subgroup of order $p^2$. Then $H'$ contains an element of order $p^2$ except in the case in which $Syl(H´)$ is isomorphic to the noncyclic subgroup of order $p^2$ of $G_p$. The corresponding subgroup of $\Hol(C_{p^2})$ is generated by $(p,Id)$ and $(0,\sigma^{p-1})$. Its intersection with $Stab_H(1)$ consists in the elements $(1-k^{\l(p-1)},\sigma^{\l(p-1)}), 1\leq \l \leq p$, since $k^{p-1} \equiv 1 \pmod{p}$, hence this intersection has order $p$. We have then that if $Syl(H´)$ is isomorphic to the noncyclic subgroup of order $p^2$ of $G_p$, then $p$ divides exactly $[H':Stab_H(1) \cap H']$ and $H'$ is not transitive. We have proved then that all transitive subgroups of $\Hol(C_{p^2})$ contain an element of order $p^2$.

 Let us look now at $Hol(C_p\times C_p)$. By \cite{Ko}, Theorem 4.4, $Hol(C_p\times C_p)$ has no elements of order $p^2$. Taking into account what we have proved above, this finishes the proof of the proposition.
 $\Box$

\begin{remark} Kohl  proves in \cite{Ko} that any Hopf Galois structure on a cyclic extension of order $p^n$, for $p$ an odd prime, is of cyclic type. Childs studies in \cite{Ch3} these Hopf Galois structures in the case of cyclic extensions of order $p^2$.
\end{remark}

\section{Description of the algorithm}

Given a separable field extension $L/K$ of degree $g$, $\wL$ its Galois closure, $G=\Gal(\wL/K),$ \newline $G'=\Gal(\wL/L)$, the action of
$G$ on the left cosets $G/G'$ is transitive, hence the morphism $\lambda:G \rightarrow S_g$ identifies $G$ with a transitive subgroup of $S_g$, which is determined up to conjugation. Moreover, if we enumerate the left cosets $G/G'$ starting with the one containing $1_G$, $\lambda(G')$ is equal to the stabilizer of $1$ in $G$. Therefore considering all separable field extensions $L/K$ of degree $g$ is equivalent to considering all transitive groups $G$ of degree $g$, up to conjugation.
The algorithm structure is as follows:

\begin{enumerate}[{Step} 1]
\item Given a transitive group $G$ of degree $g$ and a type of regular subgroups $N$ of $S_g$, run over the conjugation class of $N$ in $S_g$ and determine whether $N$ is normalized by $G$. In the affirmative case, check if $N$ is contained in $G$.
\item For each transitive group $G$ of degree $g$ and $G'=Stab_G(1)$, determine the number $intfields(G)$ of subgroups of $G$ containing $G'$, that is, by the fundamental  theorem of classical Galois theory, the number of intermediate fields of the extension $L/K$.
\item For each pair $(G,N)$ determined in Step 1, determine the number $subGst(N)$ of $G$-stable subgroups of $N$, i.e. subgroups of $N$ normalized by $G$, that is, the cardinal of the image of the map $\mathcal{F}_H$ in Theorem \ref{esto} for the Hopf Galois structure given by $N$. Check if this number equals $intfields(G)$, that is if the Galois correspondence is bijective.
\item For each pair $(G,N),(G,N')$, with $N\simeq N'$ and  $subGst(N)=subGst(N')$, check if $N$ and $N'$ are $G$-isomorphic, that is if the corresponding Hopf algebras are isomorphic. To this end, we use that for a regular subgroup $N$ of the symmetric group $S_g$, the automorphism group $\Aut N$ of $N$ is isomorphic to the stabilizer of 1 in the holomorph $Hol(N)$ of $N$ and that $Hol(N)$ is the normalizer of $N$ in $S_g$.
\end{enumerate}

We note that in Step 1 we compute the transversal of the normalizer of $N$ in $S_g$ and the conjugate of $N$ by each element in this transversal. This computation occurs to need a significantly shorter  execution time than the use of the Magma function Class from degree 9 onwards. In the vector which collects the regular subgroups $N$ of $S_g$ giving a Hopf Galois structure we have added a numbering variable. In this way, each of these $N$'s is identified with an integer number. This numeration is respected all along the program so that, once the $N$'s have been computed in Step 1, we can easily know the properties of the corresponding Hopf Galois structures by searching the assigned number. This greatly simplifies the reading and interpretation of the results.

\section{Results and discussion}

We present the results obtained for separable field extensions of degrees $g=4, 6, 8, 9$ and $10$ in the tables in appendix 1. We denote by $kTi$ the $i$th transitive group of degree $k$ called by TransitiveGroup$(k,i)$ in the Magma program. In tables \ref{d4-ext} to \ref{d10-ext}, for each transitive group $G$ of degree $g$ and each group $N$ of order $g$, the corresponding table gives the total number $T$ of Hopf Galois structures of type $N$ for a separable field extension $L/K$ of degree $g$ such that the Galois group of the normal closure $\widetilde{L}$ over $K$ is isomorphic to $G$. Moreover, it gives the number a-c of those which are almost classically Galois, the number BC of those for which the Galois correspondence is bijective and the number G-i of Hopf algebra isomorphism classes in which the Hopf Galois structures are partitioned. In particular the difference BC minus a-c gives the number of non almost classically Galois Hopf Galois structures for which the Galois correspondence is bijective. We also give as a summary the total number of Hopf Galois structures corresponding to $G$ together with the above items. The transitive groups $G$ such that the corresponding field extension $L/K$ has no Hopf Galois structure are not included in the table.

\vspace{0.3cm}
We note that the field extension with smallest degree having a non almost classically Galois Hopf Galois structure with bijective Galois correspondence is a Galois extension of degree 4 with Galois group $C_4$ and the Hopf Galois structure is of type $C_2\times C_2$. The non-Galois extension with smallest degree having this property is a separable extension of degree 6 whose Galois closure has group $6T5$ and the Hopf Galois structure is of type $S_3$. The field extension with smallest degree having non-isomorphic Hopf Galois structures with isomorphic Hopf algebras is a Galois extension of degree 6 with Galois group the symmetric group $S_3$ for which the three Hopf Galois structures of cyclic type $C_6$ have underlying isomorphic Hopf algebras.

\vspace{0.3cm}
In table \ref{d8-ext-Giso} we give the distribution of Hopf Galois structures in Hopf algebra isomorphism classes for transitive groups of degree 8 having some class with more that one element. For example, in the cell corresponding to $G=N=C_4\times C_2, 10=5\times 1+1 \times 2+1 \times 3$ means that for a Galois extension with Galois group $C_4\times C_2$ there are 10 Hopf Galois structures of type $C_4\times C_2$ which are distributed in 5 classes with 1 element, 1 class with 2 elements and 1 class with 3 elements.

\vspace{0.3cm}
Table \ref{fig} reflects the size of the results obtained and the computation complexity of the algorithm. In it we give for every degree $g$ the order of the symmetric group of degree $g$; the total number of transitive groups of degree $g$ and the number Max of transitive groups of degree $g$ whose order does not exceed the order of the holomorphs of all the groups of order $g$; the number of possible types of Hopf Galois structures; the total number of Hopf Galois structures and the number of the almost classically Galois ones; the number of Hopf Galois structures with bijective Galois correspondence and the number of those which are not almost classically Galois; the number of Hopf algebra isomorphism classes in which the Hopf Galois structures are partitioned (which correspond to $G$-isomorphism classes of the corresponding regular groups $N$) and the number of those for Galois extensions; and finally the execution times in seconds and the memory used in megabytes. We note that the presented program is very efficient up to degree 11. One may observe in particular that the computation for degree 8, which gives a large number of Hopf Galois structures, takes only about 17 seconds. For degree 12, even on a supercomputer, the computation cannot be carried out due to lack of memory.

\subsection{Conclusions}

The elaboration of the algorithm presented allows to determine all Hopf Galois structures of separable field extensions of a given degree up to degree 11. Such a determination has been obtained by theoretic tools only for prime degree extensions. In proposition \ref{psqua} we prove a partial result concerning Hopf Galois structures of separable field extensions of degree $p^2$, for $p$ an odd prime, which came up from the results for degree 9 extensions. In general Hopf Galois structures may not be determined by theoretic reasoning, therefore it is valuable to have the presented algorithm at disposal. We specially highlight the richness of results obtained in the degree 8 case.

\subsection{Example}

We consider a Galois extension $L/K$ with Galois group $G=C_2\times C_2\times C_2$. As given in table \ref{d8-ext-Giso}, it has 42 Hopf Galois structures of type $D_{2\cdot 4}$ partitioned in 7 Hopf algebra isomorphism classes of 6 elements each. We will examine in detail one of these classes and determine the corresponding Hopf algebra and Galois actions. We may write $L=K(\alpha,\beta,\gamma)$, with $\alpha^2, \beta^2, \gamma^2 \in K$ and $G$ is then generated by the automorphisms $a,b,c$ given by

$$\begin{array}{cccr} a:& \alpha & \mapsto & -\alpha \\ & \beta & \mapsto & \beta \\ &\gamma & \mapsto & \gamma \end{array}, \quad
\begin{array}{cccr} b:& \alpha & \mapsto & \alpha \\ & \beta & \mapsto & -\beta \\ &\gamma & \mapsto & \gamma \end{array}, \quad
\begin{array}{cccr} c:& \alpha & \mapsto & \alpha \\ & \beta & \mapsto & \beta \\ &\gamma & \mapsto & -\gamma \end{array}.$$

\noindent The group $C_2\times C_2\times C_2\simeq 8T3$ is given in Magma as the subgroup of the symmetric group generated by $(1,8)(2,3)(4,5)(6,7),(1,3)(2,8)(4,6)(5,7),(1,5)(2,6)(3,7)(4,8)$. If we order the elements in $G$ as $Id, ab, b, ac, c, abc, bc, a$, we have $\lambda(a)=(1,8)(2,3)(4,5)(6,7), \linebreak \lambda(b)=(1,3)(2,8)(4,6)(5,7), \lambda(c)=(1,5)(2,6)(3,7)(4,8)$ and we shall identify $G$ with its image by $\lambda$. The following regular subgroups of $S_8$ are isomorphic to $D_{2\cdot 4}$, normalized by $G$ and mutually $G$-isomorphic.

$$\begin{array}{c} N_1=\langle s_1=(1,8)(2,3)(4,5)(6,7), r_1=(1,6,5,2),(3,4,7,8) \rangle, \\ N_2=\langle s_2=(1,8)(2,3)(4,5)(6,7), r_2=(1,4,7,2),(3,6,5,8) \rangle, \\ N_3=\langle s_3=(1,6)(2,5)(3,4)(7,8), r_3=(1,4,5,8),(2,3,6,7) \rangle, \\ N_4=\langle s_4=(1,2)(3,8)(4,7)(5,6), r_4=(1,4,3,6),(2,5,8,7) \rangle, \\ N_5=\langle s_5=(1,6)(2,5)(3,4)(7,8), r_5=(1,2,3,8),(4,5,6,7) \rangle, \\ N_6=\langle s_6=(1,4)(2,7)(3,6)(5,8), r_6=(1,6,7,8),(2,3,4,5) \rangle. \end{array}$$

We check that $ar_ia=r_i^3, br_ib=r_i, cr_ic=r_i, as_ia=s_i, bs_ib=s_i, cs_ic=s_i, 1\leq i \leq 6$, hence $N_i$ is normalized by $G$, for $1\leq i \leq 6$ and $s_i \mapsto s_j, r_i \mapsto r_j$ defines a $G$-isomorphism from $N_i$ to $N_j, 1\leq i,j \leq 6$.

By computation, we obtain that the Hopf algebra corresponding to $N_i$ is the $K$-Hopf algebra with basis $1, r_i+r_i^3, r_i^2, \alpha (r_i-r_i^3), s_i, s_ir_i+s_ir_i^3, s_ir_i^2, \alpha (s_ir_i-s_ir_i^3)$ and the Hopf actions are given by

$$\begin{array}{lllllll} r_1 \mapsto ab, & r_1^2 \mapsto c, & r_1^3 \mapsto abc, & s_1 \mapsto a, & s_1r_1 \mapsto bc, & s_1r_1^2 \mapsto ac, & s_1r_1^3 \mapsto b \\
r_2 \mapsto ab, & r_2^2 \mapsto bc, & r_2^3 \mapsto ac, & s_2 \mapsto a, & s_2r_2 \mapsto c, & s_2r_2^2 \mapsto abc, & s_2r_2^3 \mapsto b \\
 r_3 \mapsto a, & r_3^2 \mapsto c, & r_3^3 \mapsto ac, & s_3 \mapsto abc, & s_3r_3 \mapsto b, & s_3r_3^2 \mapsto ab, & s_3r_3^3 \mapsto bc \\
 r_4 \mapsto abc, & r_4^2 \mapsto b, & r_4^3 \mapsto ac, & s_4 \mapsto ab, & s_4r_4 \mapsto bc, & s_4r_4^2 \mapsto a, & s_4r_4^3 \mapsto c \\
r_5 \mapsto ab, & r_5^2 \mapsto b, & r_5^3 \mapsto a, & s_5 \mapsto abc, & s_5r_5 \mapsto c, & s_5r_5^2 \mapsto ac, & s_5r_5^3 \mapsto bc \\
 r_6 \mapsto abc, & r_6^2 \mapsto bc, & r_6^3 \mapsto a, & s_6 \mapsto ac, & s_6r_6 \mapsto b, & s_6r_6^2 \mapsto ab, & s_6r_6^3 \mapsto c
\end{array}
$$

A different explicit example can be found in \cite{S}, Example 5.3.1.

\section*{Acknowledgments} We are grateful to Anna Rio and Montserrat Vela for valuable discussions on the subject of this paper and to Joan Nualart and Pawe\l \ Bogdan for their help with the Magma program. We thank also Jordi Gu\`ardia for giving us access to the supercomputer of the Universitat Polit\`ecnica de Catalunya.

\newpage

\section*{Appendix 1 - Tables}\label{tables}

\begin{table}[ht]
\caption{Degree 4 extensions}
\label{d4-ext}
\begin{center}
\begin{tabular}{|c||c|c|c|c||c|c|c|c||c|c|c|c|}
\hline
\multicolumn{1}{|c||}{} & \multicolumn{12}{|c|}{\bf Hopf Galois structures} \\
\hline \hline
 \multicolumn{1}{|c||}{\bf Galois } & \multicolumn{4}{|c||}{Type $C_4$} & \multicolumn{4}{|c||}{Type $C_2\times C_2$} & \multicolumn{4}{|c|}{Summary} \\ \cline{2-13}  \multicolumn{1}{|c||}{\bf group } & T & a-c &  BC & G-i & T & a-c  & BC & G-i & T & a-c  & BC & G-i \\
\hline $4T1 \simeq C_4$ & 1&1&1&1 & 1&0&1&1&2&1&2&2\\ \hline
$4T2 \simeq C_2\times C_2$ & 3&0&0&3 & 1&1&1&1&4&1&1&4 \\ \hline
$4T3  \simeq D_{2\cdot 4}$ & 1&1&1&1 & 1&1&1&1&2&2&2&2\\ \hline
$4T4 \simeq A_4$ & 0&0&0&0 & 1&1&1&1&1&1&1&1 \\ \hline
$4T5 \simeq S_4$ & 0&0&0&0 & 1&1&1&1&1&1&1&1 \\ \hline
\end{tabular}
\end{center}
\end{table}
\begin{table}[hb]
\caption{Degree 6 extensions}
\label{d6-ext}
\begin{center}
\begin{tabular}{|c||c|c|c|c||c|c|c|c||c|c|c|c|}
\hline
\multicolumn{1}{|c||}{} & \multicolumn{12}{|c|}{\bf Hopf Galois structures} \\
\hline \hline
 \multicolumn{1}{|c||}{\bf Galois } & \multicolumn{4}{|c||}{Type $C_6\simeq C_2\times C_3$} & \multicolumn{4}{|c||}{Type $D_{2\cdot 3}\simeq S_3$} & \multicolumn{4}{|c|}{Summary} \\ \cline{2-13}  \multicolumn{1}{|c||}{\bf group } & T & a-c &  BC & G-i & T & a-c  & BC & G-i & T & a-c  & BC & G-i \\
\hline $6T1 \simeq C_6$ & 1&1&1&1 & 2&0&1&2&3&1&2&3 \\ \hline
$6T2 \simeq D_{2\cdot 3}$ & 3&0&0&1 & 2&1&1&2&5&1&1&3 \\ \hline
$6T3  \simeq C_6\rtimes C_2$ & 1&1&1&1 & 2&1&1&2&3&2&2&3 \\ \hline
$6T5 \simeq (C_3\times C_3)\rtimes C_2$ & 0&0&0&0 & 2&1&2&2&2&1&2&2 \\ \hline
$6T9 \simeq 6T5\rtimes C_2$ & 0&0&0&0 & 2&2&2&2&2&2&2&2 \\ \hline
\end{tabular}
\end{center}
\end{table}

\begin{landscape}

\begin{table}[h]
\centering
\caption{Degree 8 extensions}
\label{d8-ext}
\vspace{0.2 cm}
\begin{tabular}{|c||c|c|c|c||c|c|c|c||c|c|c|c||c|c|c|c||c|c|c|c|}
\hline
\multicolumn{1}{|c||}{} & \multicolumn{20}{|c|}{\bf Hopf Galois structures} \\
\hline \hline
 \multicolumn{1}{|c||}{\bf Galois } & \multicolumn{4}{|c||}{Type $C_8$} & \multicolumn{4}{|c||}{Type $C_4\times C_2$} & \multicolumn{4}{|c||}{Type $C_2\times C_2\times C_2$} & \multicolumn{4}{|c||}{Type $D_{2\cdot4}$} & \multicolumn{4}{|c|}{Type $Q_8$}\\ \cline{2-21}  \multicolumn{1}{|c||}{\bf group } & T & a-c &  BC & G-i & T & a-c  & BC & G-i & T & a-c  & BC & G-i & T & a-c  & BC & G-i & T & a-c  & BC & G-i \\
\hline
\multicolumn{1}{|c||}{$8T1\simeq C_8$}                     & 2 &1&2&2& 0 &0&0&0     & 0 &0&0&0                & 2 &0&2&2     & 2 &0&2&2    \\ \hline
\multicolumn{1}{|c||}{$8T2\simeq C_4\times C_2$}           & 4 &0&0&2 & 10  &1&1&7 & 4 &0&1&4                 & 6&0&2&5    & 2  &0&0&2   \\ \hline
\multicolumn{1}{|c||}{$8T3\simeq (C_2)^3$} & 0&0&0&0  & 42 &0&0&28    & 8 &1&1&8              & 42 &0&0&7   & 14 &0&0&7   \\ \hline
\multicolumn{1}{|c||}{$8T4\simeq D_{2\cdot4}$}             & 2 &0&0&1 & 14 &0&0&9     & 6  &0&0&4             & 6 &1&1&4   & 2 &0&0&2    \\ \hline
\multicolumn{1}{|c||}{$8T5\simeq Q_8$}                     & 6&0&0&3  & 6  &0&6&3     & 2 &0&2&1              & 6 &0&6&6   & 2 &1&2&2  \\ \hline
\multicolumn{1}{|c||}{$8T6$}                    & 2 &1&2&2   & 0  &0&0&0             & 0 &0&0&0               & 2 &1&2&2   & 2 &0&2&2    \\ \hline
\multicolumn{1}{|c||}{$8T7$}                     & 2&2&2&2     & 0 &0&0&0              & 0 &0&0&0             & 2&0&2&2    & 2&0&2&2     \\ \hline
\multicolumn{1}{|c||}{$8T8$}                     & 2 &1&2&2    & 0 &0&0&0              & 0 &0&0&0             & 2&0&2&2    & 2 &1&2&2    \\ \hline
\multicolumn{1}{|c||}{$8T9$}                     & 0 &0&0&0    & 10&1&1&9           & 4&1&1&4                 & 6 &2&2&5   & 2 &0&0&2    \\ \hline
\multicolumn{1}{|c||}{$8T10$}                     & 0 &0&0&0    & 6 &2&3&6          & 4 &0&1&4       & 0 &0&0&0            & 0 &0&0&0    \\ \hline
\multicolumn{1}{|c||}{$8T11$}                     & 2 &0&0&1   & 6  &2&6&5          & 2 &0&2&2       & 6  &1&6&6           & 2 &1&2&2  \\ \hline
\multicolumn{1}{|c||}{$8T12$}                     & 0&0&0&0     & 0&0&0&0               & 2 &0&2&1   & 0 &0&0&0           & 2  &1&2&2   \\ \hline
\multicolumn{1}{|c||}{$8T13$}                     & 0 &0&0&0    & 0  &0&0&0             & 2 &1&1&2  & 0 &0&0&0            & 2 &0&0&1    \\ \hline
\multicolumn{1}{|c||}{$8T14$}                     & 0 &0&0&0   & 0&0&0&0              & 4 &0&1&3       & 0&0&0&0            & 0 &0&0&0    \\ \hline
\multicolumn{1}{|c||}{$8T15$}                     & 2 &2&2&2    & 0&0&0&0            & 0&0&0&0         & 2 &1&2&2         & 2&1&2&2     \\ \hline
\multicolumn{1}{|c||}{$8T16$}                     & 0&0&0&0    & 0 &0&0&0              & 0 &0&0&0      & 2 &0&2&2         & 2&0&2&2     \\ \hline
\multicolumn{1}{|c||}{$8T17$}                     & 0&0&0&0     & 0&0&0&0               & 0 &0&0&0     & 2  &1&2&2           & 2&1&2&2     \\ \hline
\multicolumn{1}{|c||}{$8T18$}                     & 0&0&0&0     & 6 &3&3&6              & 4 &1&1&4     & 0 &0&0&0            & 0 &0&0&0    \\ \hline

\end{tabular}
\end{table}

\begin{table}[h]
\centering
\caption{Degree 8 extensions (cont.)}
\label{d8-ext-cont}
\vspace{0.2 cm}
\begin{tabular}{|c||c|c|c|c||c|c|c|c||c|c|c|c||c|c|c|c||c|c|c|c|}
\hline
\multicolumn{1}{|c||}{} & \multicolumn{20}{|c|}{\bf Hopf Galois structures} \\
\hline \hline
 \multicolumn{1}{|c||}{\bf \quad \, Galois \quad \, \, } & \multicolumn{4}{|c||}{Type $C_8$} & \multicolumn{4}{|c||}{Type $C_4\times C_2$} & \multicolumn{4}{|c||}{Type $C_2\times C_2\times C_2$} & \multicolumn{4}{|c||}{Type $D_{2\cdot4}$} & \multicolumn{4}{|c|}{Type $Q_8$}\\ \cline{2-21}  \multicolumn{1}{|c||}{\bf group } & T & a-c &  BC & G-i & T & a-c  & BC & G-i & T & a-c  & BC & G-i & T & a-c  & BC & G-i & T & a-c  & BC & G-i \\
\hline
\multicolumn{1}{|c||}{ $8T19$}                     & 0 &0&0&0    & 2 &1&2&2              & 2 &1&2&2     & 0&0&0&0            & 0 &0&0&0    \\ \hline
\multicolumn{1}{|c||}{$8T20$}                     & 0&0&0&0     & 2 &0&2&2              & 2 &0&2&2     & 0 &0&0&0            & 0 &0&0&0    \\ \hline
\multicolumn{1}{|c||}{$8T22$}                     & 0&0&0&0     & 6 &6&6&6             &2 &2&2&2       & 6  &6&6&6            & 2  &2&2&2  \\ \hline
\multicolumn{1}{|c||}{$8T23$}                     & 0&0&0&0     & 0 &0&0&0              & 0 &0&0&0     & 0  &0&0&0          & 2 &1&2&2     \\ \hline
\multicolumn{1}{|c||}{$8T24$}                     & 0&0&0&0    & 0&0&0&0               & 2 &1&1&2     & 0 &0&0&0           & 0 &0&0&0    \\ \hline
\multicolumn{1}{|c||}{$8T25$}                     & 0  &0&0&0   & 0  &0&0&0             & 1  &1&1&1    & 0 &0&0&0           & 0 &0&0&0   \\ \hline
\multicolumn{1}{|c||}{$8T26$}                     & 0 &0&0&0    & 0 &0&0&0              & 0  &0&0&0     & 2  &2&2&2           & 2  &2&2&2 \\ \hline
\multicolumn{1}{|c||}{$8T29$}                     & 0 &0&0&0    & 2 &2&2&2              & 2  &2&2&2      & 0 &0&0&0           & 0 &0&0&0    \\ \hline
\multicolumn{1}{|c||}{$8T32$}                     & 0  &0&0&0   & 0   &0&0&0            & 2  &2&2&2      & 0 &0&0&0           & 2  &2&2&2   \\ \hline
\multicolumn{1}{|c||}{$8T33$}                     & 0 &0&0&0    & 0 &0&0&0              & 1 &1&1&1       & 0 &0&0&0           & 0  &0&0&0   \\ \hline
\multicolumn{1}{|c||}{$8T34$}                     & 0 &0&0&0    & 0 &0&0&0              & 3 &0&3&3       & 0 &0&0&0           & 0 &0&0&0    \\ \hline
\multicolumn{1}{|c||}{$8T36$}                     & 0 &0&0&0    & 0  &0&0&0             & 1 &1&1&1       & 0  &0&0&0          & 0  &0&0&0   \\ \hline
\multicolumn{1}{|c||}{$8T37$}                     & 0&0&0&0     & 0 &0&0&0              & 2 &0&2&2      & 0 &0&0&0           & 0&0&0&0     \\ \hline
\multicolumn{1}{|c||}{$8T39$}                     & 0 &0&0&0    & 0 &0&0&0              & 2&2&2&2                         & 0  &0&0&0          & 0&0&0&0     \\ \hline
\multicolumn{1}{|c||}{$8T40$}                     & 0 &0&0&0    & 0 &0&0&0              & 0 &0&0&0        & 0&0&0&0            & 2&2&2&2     \\ \hline
\multicolumn{1}{|c||}{$8T41$}                     & 0 &0&0&0    & 0  &0&0&0             & 1 &1&1&1                        & 0 &0&0&0           & 0 &0&0&0    \\ \hline
\multicolumn{1}{|c||}{$8T48$}                     & 0 &0&0&0    & 0  &0&0&0             & 1 &1&1&1      & 0 &0&0&0           & 0 &0&0&0    \\ \hline
\end{tabular}
\end{table}

\begin{table}[h]
\centering
\caption{Degree 8 extensions-Summary}
\label{d8-ext-sum}
\vspace{0.2 cm}
\begin{tabular}{|c||c|c|c|c|}
\hline
\multicolumn{1}{|c||}{} & \multicolumn{4}{|c|}{\bf H G structures} \\
\hline \hline
 \multicolumn{1}{|c||}{\bf Galois group} &  T & a-c  & BC & G-i \\
\hline
\multicolumn{1}{|c||}{$8T1\simeq C_8$}                     & 6 &1&6&6   \\ \hline
\multicolumn{1}{|c||}{$8T2\simeq C_4\times C_2$}           & 26 &1&4&20  \\ \hline
\multicolumn{1}{|c||}{$8T3\simeq (C_2)^3$} & 106&1&1&50  \\ \hline
\multicolumn{1}{|c||}{$8T4\simeq D_{2\cdot4}$}  &30 &1&1&20   \\ \hline
\multicolumn{1}{|c||}{$8T5\simeq Q_8$}                     & 22&1&16&15  \\ \hline
\multicolumn{1}{|c||}{$8T6$}                    & 6 &2&6&6  \\ \hline
\multicolumn{1}{|c||}{$8T7$}                     &  6 &2&6&6     \\ \hline
\multicolumn{1}{|c||}{$8T8$}                     &  6 &2&6&6   \\ \hline
\multicolumn{1}{|c||}{$8T9$}                     & 22 &4&4&20    \\ \hline
\multicolumn{1}{|c||}{$8T10$}                     & 10 &2&4&10   \\ \hline
\multicolumn{1}{|c||}{$8T11$}                     & 18 &4&16&16  \\ \hline
\multicolumn{1}{|c||}{$8T12$}                     & 4&1&4&3   \\ \hline
\multicolumn{1}{|c||}{$8T13$}                     & 4 &1&1&3    \\ \hline
\multicolumn{1}{|c||}{$8T14$}                     & 4 &0&1&3    \\ \hline
\multicolumn{1}{|c||}{$8T15$}                     & 6 &4&6&6     \\ \hline
\multicolumn{1}{|c||}{$8T16$}                     & 4&0&4&4    \\ \hline
\multicolumn{1}{|c||}{$8T17$}                     & 4&2&4&4     \\ \hline
\multicolumn{1}{|c||}{$8T18$}                     & 10&4&4&10    \\ \hline

\end{tabular}
\hspace{1cm}
\begin{tabular}{|c||c|c|c|c|}
\hline
\multicolumn{1}{|c||}{} & \multicolumn{4}{|c|}{\bf H G structures} \\
\hline \hline
 \multicolumn{1}{|c||}{\bf Galois group} &  T & a-c &  BC & G-i \\
\hline
\multicolumn{1}{|c||}{ $8T19$}                     & 4 &2&4&4  \\ \hline
\multicolumn{1}{|c||}{$8T20$}                     & 4&0&4&4 \\ \hline
\multicolumn{1}{|c||}{$8T22$}                     & 16&16&16&16  \\ \hline
\multicolumn{1}{|c||}{$8T23$}                     & 2&1&2&2    \\ \hline
\multicolumn{1}{|c||}{$8T24$}                     & 2&1&1&2       \\ \hline
\multicolumn{1}{|c||}{$8T25$}                     & 1  &1&1&1   \\ \hline
\multicolumn{1}{|c||}{$8T26$}                     & 4 &4&4&4 \\ \hline
\multicolumn{1}{|c||}{$8T29$}                     & 4 &4&4&4  \\ \hline
\multicolumn{1}{|c||}{$8T32$}                     & 4 &4&4&4 \\ \hline
\multicolumn{1}{|c||}{$8T33$}                     &  1 &1&1&1   \\ \hline
\multicolumn{1}{|c||}{$8T34$}                     & 3 &0&3&3    \\ \hline
\multicolumn{1}{|c||}{$8T36$}                     & 1 &1&1&1   \\ \hline
\multicolumn{1}{|c||}{$8T37$}                     &  2 &0&2&2   \\ \hline
\multicolumn{1}{|c||}{$8T39$}                     &  2&2&2&2      \\ \hline
\multicolumn{1}{|c||}{$8T40$}                     & 2&2&2&2     \\ \hline
\multicolumn{1}{|c||}{$8T41$}                     &1 &1&1&1   \\ \hline
\multicolumn{1}{|c||}{$8T48$}                     & 1 &1&1&1  \\ \hline
\multicolumn{5}{c}{}
\end{tabular}
\end{table}

\begin{table}[h]
\centering
\caption{Hopf algebra isomorphism classes for degree 8 extensions}
\label{d8-ext-Giso}
\vspace{0.2 cm}
\begin{tabular}{|c||c|c|c|c|c|}
\hline
\multicolumn{1}{|c||}{} & \multicolumn{5}{|c|}{\bf Isomorphism classes} \\
\hline \hline
 \multicolumn{1}{|c||}{\bf{Galois group}} & \multicolumn{1}{|c|}{Type $C_8$} & \multicolumn{1}{|c|}{Type $C_4\times C_2$} & \multicolumn{1}{|c|}{Type $C_2\times C_2\times C_2$} & \multicolumn{1}{|c|}{Type $D_{2\cdot4}$} & \multicolumn{1}{|c|}{Type $Q_8$}\\ \hline
\multicolumn{1}{|c||}{$8T2\simeq C_4\times C_2$}           & $4=2\times 2$ & $10=5\times 1+1 \times 2+1 \times 3$ & $4=4\times 1$&$6=4\times 1+1 \times2$& $2=2\times 1$  \\ \hline
\multicolumn{1}{|c||}{$8T3\simeq (C_2)^3$} & 0  & $42=21\times 1+7\times 3$    & $8=8\times 1$ &$42=7\times 6$ & $14=7\times 2$   \\ \hline
\multicolumn{1}{|c||}{$8T4\simeq D_{2\cdot4}$}  & $2=1\times 2$ & $14=4\times 1+5\times 2$ & $6=2\times 1+2\times 2$ & $6=3\times 1+1\times 3$ & $2=2\times 1$   \\ \hline
\multicolumn{1}{|c||}{$8T5\simeq Q_8$} &$6=3\times 2$ &$6=3\times 2$ & $2=1\times 2$ &$6=6\times 1$ &$2=2\times 1$ \\ \hline
\multicolumn{1}{|c||}{$8T9$} &$0$&$10=8\times 1+1\times 2$&$4=4\times 1$&$6=4\times 1+1\times 2$&$2=2\times 1$    \\ \hline
\multicolumn{1}{|c||}{$8T11$}  &$2=1\times 2$&$6=4\times 1+1\times 2$&$2=2\times 1$&$6=6\times 1$&$2=2\times 1$    \\ \hline
\multicolumn{1}{|c||}{$8T12$} &$0$&$0$&$2=1\times 2$&$0$&$2=2\times 1$    \\ \hline
\multicolumn{1}{|c||}{$8T13$}   & $0$ &$0$&$2=2\times 1$&$0$&$2=1\times 2$    \\ \hline
\multicolumn{1}{|c||}{$8T14$}                     & $0$ &$0$&$4=2\times 1+1\times 2$&$0$    & $0$    \\ \hline
\end{tabular}
\end{table}

\begin{table}[ht]
\caption{Degree 9 extensions}
\label{d9-ext}
\begin{center}
\begin{tabular}{|c||c|c|c|c||c|c|c|c||c|c|c|c|}
\hline
\multicolumn{1}{|c||}{} & \multicolumn{12}{|c|}{\bf Hopf Galois structures} \\
\hline
\hline
 \multicolumn{1}{|c||}{\bf Galois } & \multicolumn{4}{|c||}{Type $C_9$} & \multicolumn{4}{|c|}{Type $C_3\times C_3$}  & \multicolumn{4}{|c|}{Summary}\\ \cline{2-13}  \multicolumn{1}{|c||}{\bf group } & T & a-c &  BC & G-i & T & a-c  & BC & G-i& T & a-c  & BC & G-i \\
\hline $9T1 \simeq C_9$ & 3&1&3&3 & 0&0&0&0& 3&1&3&3 \\ \hline
$9T2 \simeq C_3\times C_3$ & 0&0&0&0 & 9&1&1&5 & 9&1&1&5\\ \hline
$9T3$ & 1&1&1&1 & 0&0&0&0 & 1&1&1&1\\ \hline
$9T4$ & 0&0&0&0 & 3&1&1&2 & 3&1&1&2\\ \hline
$9T5$ & 0&0&0&0 & 1&1&1&1& 1&1&1&1 \\ \hline
$9T6$ & 3&3&3&3 & 0&0&0&0 & 3&3&3&3\\ \hline
$9T7$ & 0&0&0&0 & 3&3&3&3 & 3&3&3&3\\ \hline
$9T8$ & 0&0&0&0 & 1&1&1&1 & 1&1&1&1\\ \hline
$9T9$ & 0&0&0&0 & 1&1&1&1 & 1&1&1&1\\ \hline
$9T10$ & 1&1&1&1& 0&0&0&0 & 1&1&1&1\\ \hline
$9T11$ & 0&0&0&0 & 1&1&1&1 & 1&1&1&1\\ \hline
$9T12$ & 0&0&0&0 & 3&3&3&3 & 3&3&3&3\\ \hline
$9T13$ & 0&0&0&0 & 1&1&1&1 & 1&1&1&1\\ \hline
$9T14$ & 0&0&0&0 & 1&1&1&1 & 1&1&1&1\\ \hline
$9T15$ & 0&0&0&0 & 1&1&1&1 & 1&1&1&1\\ \hline
$9T16$ & 0&0&0&0 & 1&1&1&1& 1&1&1&1 \\ \hline
$9T18$ & 0&0&0&0 & 1&1&1&1 & 1&1&1&1\\ \hline
$9T19$ & 0&0&0&0 & 1&1&1&1 & 1&1&1&1\\ \hline
$9T23$ & 0&0&0&0 & 1&1&1&1 & 1&1&1&1\\ \hline
$9T26$ & 0&0&0&0 & 1&1&1&1 & 1&1&1&1\\ \hline
\end{tabular}
\end{center}
\end{table}

\begin{table}[ht]
\caption{Degree 10 extensions}
\label{d10-ext}
\begin{center}
\begin{tabular}{|c||c|c|c|c||c|c|c|c||c|c|c|c|}
\hline
\multicolumn{1}{|c||}{} & \multicolumn{12}{|c|}{\bf Hopf Galois structures} \\
\hline
\hline
 \multicolumn{1}{|c||}{\bf Galois } & \multicolumn{4}{|c||}{Type $C_{10}$} & \multicolumn{4}{|c|}{Type $D_{2\cdot 5}$}  & \multicolumn{4}{|c|}{Summary}\\ \cline{2-13}  \multicolumn{1}{|c||}{\bf group } & T & a-c &  BC & G-i & T & a-c  & BC & G-i& T & a-c  & BC & G-i \\
\hline $10T1 \simeq C_{10}$ & 1&1&1&1 & 2&0&1&2& 3&1&2&3 \\ \hline
$10T2 \simeq D_{2\cdot 5}$ & 5&0&0&1 & 2&1&1&2 & 7&1&1&3\\ \hline
$10T3$ & 1&1&1&1 & 2&1&1&2 & 3&2&2&3\\ \hline
$10T4$ & 1&0&1&1 & 2&0&1&2 & 3&0&2&3\\ \hline
$10T5$ & 1&1&1&1 & 2&1&1&2& 3&2&2&3 \\ \hline
$10T6$ & 0&0&0&0 & 2&1&2&2 & 2&1&2&2\\ \hline
$10T9$ & 0&0&0&0 & 2&2&2&2 & 2&2&2&2\\ \hline
$10T10$ & 0&0&0&0 & 2&0&2&2 & 2&0&2&2\\ \hline
$10T17$ & 0&0&0&0 & 2&2&2&2 & 2&2&2&2\\ \hline
\end{tabular}
\end{center}
\end{table}

\begin{table}[ht]
\centering
\caption{Computation figures}
\label{fig}
\vspace{0.2 cm}
\begin{tabular}{|c||c||c|c||c||c|c||c|c||c|c||c||c|}
\hline
\multicolumn{1}{|c||}{\bf Degree} & \multicolumn{1}{|c||}{\bf Order}& \multicolumn{2}{|c||}{\bf Transitive Groups} &\multicolumn{1}{|c||}{\bf Types} &\multicolumn{2}{|c||}{\bf HG struct.} &\multicolumn{2}{|c||}{\bf BC} &\multicolumn{2}{|c||}{\bf $G$-iso} & \multicolumn{1}{|c||}{\bf Execution time} &\multicolumn{1}{|c|}{\bf Memory used} \\
\cline{3-4} \cline{6-11}
 \multicolumn{1}{|c||}{} &  \multicolumn{1}{|c||}{}& \multicolumn{1}{|c|}{\quad Total\quad} & \multicolumn{1}{|c||}{Max} & \multicolumn{1}{|c||}{} & \multicolumn{1}{|c|}{Total} & \multicolumn{1}{|c||}{a-c}& \multicolumn{1}{|c|}{Total} & \multicolumn{1}{|c||}{not a-c}& \multicolumn{1}{|c|}{Total}& \multicolumn{1}{|c||}{Galois}&\multicolumn{1}{|c||}{(s)}&\multicolumn{1}{|c|}{(MB)}\\ \hline \hline
 2&2&1&1&1&1&1&1&0&1&1&$\approx 1$ & $\approx 10$  \\ \hline
3&6&2&2&1&2&2&2&0&2&1&$\approx 1$ & $\approx 10$  \\ \hline
4&24&5&5&2&10&6&7&1&10&6&$\approx 1$ & $\approx 11$ \\ \hline
5&120&5&3&1& 3&3&3&0&3&1&$\approx 1$ &$\approx 11$ \\ \hline
6& 720&16&10&2&15&7&9&2&13&6&$\approx 2$ & $\approx 11$\\ \hline
7& 5,040&7&4&1&4&4&4&0&4&1&$\approx 1$ &$\approx 11$ \\ \hline
8& 40,320&50&48& 5&348&74&147&73&262&111&$\approx 17$&$\approx 40$ \\ \hline
9& 362,880&34&26&2&38&26&28&2&33&8&$\approx 10$& $\approx 16$ \\ \hline
10&3,628,800&45&21&2&27&11&17&6&23&6&$\approx 160$&$\approx 45$ \\ \hline
11& 39,916,800&8&4&1&4&4&4&0&4&1&$\approx 90$ & $\approx 160$ \\ \hline
\end{tabular}
\end{table}

\end{landscape}

\section*{Appendix 2 - Magma code}
\begin{verbatim}
//MAIN FUNCTION
HopfGalois:=function(g)

S:=Sym(g); n:=NumberOfTransitiveGroups(g);

//We calculate triv (number of regular subgroups) and count_max
//(biggest transitive group which may have a Hopf Galois structure)
i:=g; max:=Order(S);
m:=1; ord_hol:=[]; count_max:=0;
while i le max do
  count:=0;
  while (m le n) and (Order(TransitiveGroup(g,m)) eq i) do
    m:=m+1;
    count:=count+1;
  end while;
  if count ne 0 then
    if i eq g then
      triv:=count; //triv means that L|K is Galois
      for j in [1..triv] do
        TG:=TransitiveGroup(g,j);
        Append(~ord_hol,Order(Holomorph(TG)));
      end for;
      //We calculate the maximum of ord_hol
      max:=ord_hol[1];
      for j in [2..#ord_hol] do
        //#ord_hol le triv, but it will be strictly minor
        //when there are repeated values in ord_hol
        if max lt ord_hol[j] then
          max:=ord_hol[j];
        end if;
      end for;
    end if;
    if i le max then
      count_max:=count_max+count;
    end if;
    delete(count);
  end if;
  i:=i+g;
end while;
delete(m); n:=count_max;

//-----------------------------------------------
//----DETERMINATION OF HOPF GALOIS STRUCTURES----
//-----------------------------------------------
TG:=[TransitiveGroup(g,i) : i in [1..triv]]; //Reg. subg. of Sg
NTG:=[Normalizer(S,TG[i]) : i in [1..triv]];
T:=[Transversal(S,NTG[i]) : i in [1..triv]];
Trans:=[[x : x in T[i]]   : i in [1..triv]];
//Transversal calculates right cosets; since we want the left ones,
//we calculate TGij doing (Trans[i][j]^-1)*x*Trans[i][j]
//instead of Trans[i][j]*x*(Trans[i][j]^-1)
set_i:=[i: i in [1..triv]];
useful:=[]; //Transitive groups which do have some H-G structure
HG_str:=[];     //Information about H-G structures
enes:=[]; m:=0; //We enumerate H-G structures
total:=0;       //We count the total number of H-G structures
total_ac:=0;    //We count the total number of a-c structures

for k in [1..n] do
  G:=TransitiveGroup(g,k); //it is contained in Sg
  GenG:=[x:x in Generators(G)];
  flag:=0; //it indicates if G has some HG structure
  //REGULAR SUBGROUPS OF Sg NORMALIZED BY G
  i_0:=1; aux:=[];
  while (i_0 le #set_i) and (#set_i ge 1) do
    i:=set_i[i_0];
    if Order(G) le ord_hol[i] then
      count:=0; //We count Hopf Galois structures
      count_ac:=0; //We count almost-classical structures
      for j in [1..#T[i]] do //#T[i]=#conjugacy classes of TG[i]
        GenN:=[(Trans[i][j]^-1)*x*Trans[i][j]:x in Generators(TG[i])];
        N:=sub<S|GenN>; //candidates to give a Hopf Galois structure
        //we see whether N is normalized by G
        indicator:=1; u:=1;
        while (u le #GenG) and (indicator eq 1) do
          v:=1;
          while (v le #GenN) and (indicator eq 1) do
            if GenN[v]^GenG[u] notin N then //n^g=g^-1*n*g
              indicator:=0;
            end if;
            v:=v+1;
          end while; //v
          u:=u+1; delete(v);
        end while; //u
        delete(u);
        if indicator eq 1 then
          count:=count+1;
          m:=m+1; Append(~enes,<m,k,i,N,0>);
          //PARTICULAR CASE: ALMOST-CLASSICAL GALOIS STRUCTURES
          if Centralizer(S,N) subset G then
            count_ac:=count_ac+1;
            enes[m][5]:=2;
          end if;
        end if;
        delete(N); delete(GenN); delete(indicator);
      end for; //j

      if count ne 0 then //there are H-G structures
        flag:=1; Append(~HG_str,<k,i,count,count_ac>);
      end if;
      total:=total+count; total_ac:=total_ac+count_ac;
      delete(count); delete(count_ac);
    else
      Append(~aux,i);
      //If Order(G)>Order(Hol(TG[i])) for a certain i, it will still
      //happen for the next G's, so we take away this index from set_i
    end if;
    i_0:=i_0+1;
  end while; //i
  if flag eq 1 then //there are H-G structures
    Append(~useful,k);
  end if;

  if #aux gt 0 then
    for j in aux do
      Exclude(~set_i,j);
    end for;
  end if;
  delete(G); delete(GenG); delete(flag); delete(i_0); delete(aux);
end for;//k
delete(TG); delete(NTG); delete(T); delete(Trans); delete(m); delete(set_i);

//-----------------------------------------------
//-----------BIJECTIVE CORRESPONDENCES-----------
//-----------------------------------------------
//-----CLASSICAL GALOIS THEORY-----
intfields:=[];
for k in useful do
  count:=0;
  G:=TransitiveGroup(g,k); H:=Stabilizer(G,1); //H=G'
  SC:=SubgroupClasses(G); //conjugacy classes of subgroups
  L:=[x: x in Class(G,SC[j]`subgroup), j in [1..#SC]];
  //L contains all the subgroups of G
  for j in [1..#L] do
    if H subset L[j] then
      count:=count +1;
    end if;
  end for; //j
  Append(~intfields,count);
  delete(count); delete(G); delete(H); delete(SC); delete(L);
end for; //k

//-----HOPF GALOIS THEORY-----
//We know that the Galois correspondence of the H-G theory is injective,
//but we wonder when it is bijective, as in the classical Galois theory.
subGst:=[]; //number of G-stable subgrups for each N
for m in [1..#enes] do
  if enes[m][5] eq 0 then //HG structures which are not almost-classical
    count:=0;
    k:=enes[m][2];
    G:=[x: x in Generators(TransitiveGroup(g,k))];
    N:=enes[m][4];
    SC:=SubgroupClasses(N); //conjugacy classes of subgroups
    L:=[x: x in Class(N,SC[j]`subgroup), j in [1..#SC]];
    //L contains all the subgroups of N
    for j in [1..#L] do
      //we see whether L[j] is normalized by G
      Lj:=[x: x in Generators(L[j])];
      indicator:=1; u:=1;
      while (u le #G) and (indicator eq 1) do
        v:=1;
        while (v le #Lj) and (indicator eq 1) do
          if Lj[v]^G[u] notin L[j] then //n^g=g^-1*n*g
            indicator:=0;
          end if;
          v:=v+1;
        end while; //v
        u:=u+1; delete(v);
      end while; //u
      delete(u);

      if indicator eq 1 then
        count:=count +1;
      end if;
    end for;//j
    Append(~subGst,count);
    if intfields[Position(useful,k)] eq count then
      enes[m][5]:=1;
    end if;
    delete(count); delete(G); delete(N); delete(SC); delete(L);
    delete(Lj);
    else //almost-classical structures
      k:=enes[m][2];
      Append(~subGst,intfields[Position(useful,k)]); delete(k);
  end if;
end for; //m
delete(intfields);

//-----------------------------------------------
//-------------G-ISOMORPHISM CLASSES-------------
//-----------------------------------------------
Giso:=[]; Giso_count:=[]; m:=1;
for number in [1..#HG_str] do
  k:=HG_str[number][1]; i:=HG_str[number][2]; nNiso:=HG_str[number][3];
  G:=[x:x in Generators(TransitiveGroup(g,k))];
  //Fixed k and i, nNiso gives the number of isomorph N's
  if nNiso eq 1 then
    Append(~Giso,<k,i,{<m,enes[m][4]>}>); Append(~Giso_count,<k,i,1>);
  else //nNiso ge 2
    count:=0;
    set_a:=[x: x in [m..m+nNiso-1]]; a:=set_a[1];
    while (#set_a ge 2) and (a lt set_a[#set_a]) do
      index:=[a];
      for b in [set_a[2]..set_a[#set_a]] do
        N1:=enes[a][4];
        if subGst[a] eq subGst[b] then
          N2:=enes[b][4];
          _,h:=IsIsomorphic(N1,N2); //It is always true
          H:=Stabilizer(Normalizer(S,N2),1); //Aut(N2)
          H:=[x:x in H]; //It is a subgroup of S=Sym(g)
          GN1:=[x:x in Generators(N1)];

          u:=1; indicator:=0;
          while (u le #H) and (indicator eq 0) do
            v:=1; flag:=1;
            while (v le #GN1) and (flag eq 1) do //#GN1=#Generators(N1)
              w:=1;
              while (w le #G) and (flag eq 1) do //#G=#Generators(G)
                x1:=G[w]*GN1[v]*G[w]^-1;
                //It is in N1 because N1 is normalized by G
                x2:=G[w]*H[u]*h(GN1[v])*H[u]^-1*G[w]^-1;
                if H[u]*h(x1)*(H[u]^-1) ne x2 then
                  flag:=0;
                end if;
                w:=w+1;
              end while; //w
              v:=v+1; delete(w);
            end while; //v
            delete(v);

            if flag eq 1 then
              indicator:=1; Append(~index,b);
              //a and b correspond to N1 and N2, respectively
            end if;
            u:=u+1; delete(flag);
          end while; //u
          delete(N1); delete(N2); delete(h); delete(H);
          delete(u); delete(indicator);
        end if;
      end for; //b

      Append(~Giso,<k,i,{<j,enes[j][4]>:j in index}>);
      //index keeps an isomorphism class
      count:=count+1;
      for j in index do
        Exclude(~set_a,j);
      end for;
      delete(index);
      if #set_a ne 0 then
        a:=set_a[1];
        if #set_a eq 1 then
          Append(~Giso,<k,i,{<a,enes[a][4]>}>);
          count:=count+1;
        end if;
      end if;
    end while;

    if count ne 0 then //There are G-isomorphism classes
      Append(~Giso_count,<k,i,count>);
    end if;
    delete(count);
  end if;
  m:=m+nNiso;
  delete(k); delete(i); delete(nNiso); delete(G);
end for; //j
delete(m); delete(subGst);

return <n,triv,HG_str,total,total_ac,useful,enes,Giso,Giso_count>;
end function;

//Recall that k refers to a transitive group of Sg, i -> HG struct type
//and N -> regular subgroup of Sg normalized by G

//DESCRIPTIVE FUNCTION
HGdescription:=function(g)
  HG:=HopfGalois(g);
  S:=Sym(g); n:=NumberOfTransitiveGroups(g);
  print "The number of transitive groups of", S, "is", n;
  print "From index", HG[1], "on, the transitive groups",
  "do not have any Hopf Galois structure.";
  print "-----------------------------------------";

  //Regular subgroups = subgroups G=Gal(L'|K) when L'=L
  print "The regular subgroups of S are";
  for i in [1..HG[2]] do //HG[2]=triv
    TransitiveGroup(g,i);
  end for;
  print "---------";
  print "In total, there are", HG[2], "regular subgroups.";

  //HG[3]=HG_str={x: x=<k,i,count,count_ac>}
  print "--------------------------";
  print "--------------------------";
  print "SUMMARIZED INFORMATION";
  number:=1; a:=1; b:=1;
  for k in HG[6] do //HG[6]=useful
    print "--------------------------";
    print "--------------------------";
    print "    Transitive Group", k;
    print "--------------------------";
    G:=TransitiveGroup(g,k); G;
    print "--------------------------";
    print "--------------------------";
    c:=0; c_total:=0; c_ac:=0; c_bc:=0; c_Giso:=0;

    //We count the number of types of HG structures (index i) for each k
    while (number le #HG[3]) and (HG[3][number][1] eq k) do
      number:=number+1; c:=c+1;
    end while;
    number:=number-c;
    for j in [1..c] do
      i:=HG[3][number][2];
      //The number of subgroups conjugated to N normalized by G
      //is equivalent to the number of Hopf Galois structures of type N
      print "The number of Hopf Galois structures of type";
      TransitiveGroup(g,i); print "is", HG[3][number][3];
      c_total:=c_total+HG[3][number][3];
      print "From these, the number of almost classical structures is",
      HG[3][number][4]; c_ac:=c_ac+HG[3][number][4];
      //We count the number of bijective correspondences: HG[7]=bij_cor
      count:=0;
      while a le #HG[7] and HG[7][a][2] eq k and HG[7][a][3] eq i do
        if HG[7][a][5] ne 0 then
          count:=count+1;
        end if;
        a:=a+1;
      end while;
      if count ne 0 then
        print "The number of bijective correspondences is", count;
        c_bc:=c_bc+count;
      end if;
      delete(count);
      //We recover the number of G-isomorphism classes: HG[9]=Giso_count
      print "The number of G-isomorphism classes is", HG[9][b][3];
      c_Giso:=c_Giso+HG[9][b][3]; b:=b+1;
      print "--------------------------";
      number:=number+1;
    end for; //j
    print "All in all, we conclude that for the transitive group", k;
    print "- The number of Hopf Galois structures is", c_total;
    print "  From these, the number of almost classical structures is",
    c_ac;
    print "- The number of bijective correspondences is", c_bc;
    print "- The number of G-isomorphism classes is", c_Giso;
    delete(c); delete(c_total); delete(c_ac); delete(c_bc);
    delete(c_Giso);
  end for; //k
  delete(number); delete(a); delete(b);

  //HG[4]=total, HG[5]=total_ac
  print "--------------------------";
  print "--------------------------";
  print "The total number of Hopf Galois structures is", HG[4];
  print "From these, the total number of almost classical structures is",
  HG[5];

  //HG[6]=useful
  print "-----------------------------------------";
  print "The transitive groups which do have some H-G structure are";
  HG[6];
  print "-----------------------------------------";

  //HG[7]=enes={x: x=<m,k,i,N,j>}, where j in {0,1,2}
  //0->not BC, 1->BC, 2->a-c struct (and so, BC)
  print "-----------------------------------------";
  print "----------------ENES---------------------";
  HG[7];

  print "-----------------------------------------";
  print "--------BIJECTIVE CORRESPONDENCES--------";
  print "-----------------------------------------";
  bc:=[]; //HG which are not almost-classical but do have BC
  ac:=[]; //almost-classical structures
  for m in [1..#HG[7]] do
    if HG[7][m][5] eq 1 then
      Append(~bc,HG[7][m][1]);
    elif HG[7][m][5] eq 2 then
       Append(~ac,HG[7][m][1]);
    end if;
  end for;

  print "Bijective correspondences of HG structures which are not a-c";
  bc;
  print "-----------";
  print "There are", #bc, "bijective correspondences";
  print "------------------------";
  print "Almost classical structures";
  ac;
  print "-----------";
  print "There are", #ac, "almost classical structures";
  print "------------------------";
  print "The total number of bijective Galois correspondences is",
  #bc+#ac; delete(bc); delete(ac);
  print "-----------------------------------------";

  //HG[8]=Giso={x: x=<k,i,{<m_1,Nm_1>,...,<m_r,Nm_r>}>}
  //HG[9]=Giso_count={x: x=<k,i,count>}
  print "----PARTITION IN G-ISOMORPHISM CLASSES---";
  print "-----------------------------------------";
  number:=1; y:=1; count_Galois:=0;
  print "--------------------------";
  for k in HG[6] do //HG[6]=useful
    print "--------------------------";
    print "    Transitive Group", k;
    print "--------------------------";
    print "--------------------------";
    c:=0;
    //We count the number of types of HG structures (index i) for each k
    while number le #HG[9] and HG[9][number][1] eq k do
      number:=number+1; c:=c+1;
    end while;
    number:=number-c;
    for j in [1..c] do
      i:=HG[9][number][2];
      print "Hopf Galois structures of type", i;
      count:=HG[9][number][3];
      if k le HG[2] then
        count_Galois:=count_Galois+count;
      end if;
      for p in [1..count] do
        HG[8][y][3]; y:=y+1;
      end for;
      print "The number of G-isomorphism classes is", count;
      print "--------------------------";
      number:=number+1;
      delete(i); delete(count);
    end for; //j
    delete(c);
  end for; //k
  delete(number); delete(y);

  print "-----------------------------------------";
  print "The total number of G-isomorphism classes is", #HG[8];
  print "From these, the total number corresponding to Galois ext. is",
  count_Galois; delete(count_Galois);
  return "----------------------------------------";

end function;
\end{verbatim}
\end{document}